\documentclass[11pt, reqno, psamsfonts]{amsart}
\pdfoutput=1

\usepackage{amssymb}
\usepackage{amsthm}
\usepackage{amsmath}
\usepackage{latexsym}
\usepackage[T2A,OT1]{fontenc}
\usepackage[utf8]{inputenc}
\usepackage[russian, french, english]{babel}

\usepackage{mathtools}
\usepackage[hidelinks]{hyperref}
\usepackage{amsbsy}
\usepackage{enumitem}
\usepackage{mathrsfs}
\usepackage{mathabx}
\usepackage{lmodern}
\usepackage{titlesec}
\usepackage[spacing=true,kerning=true,babel=true,tracking=true]{microtype}
\usepackage[shortcuts]{extdash}
\usepackage[foot]{amsaddr}
\usepackage[left=1in,right=1in,top=1in,bottom=1in,bindingoffset=0cm]{geometry}
\usepackage[backend=biber, style=alphabetic, sorting=nyt, maxnames=100]{biblatex}
\title{Improved Lower Bound for Difference Bases}
\date{}
\author{Anton~Bernshteyn}
\email{abernsht@math.cmu.edu}
\author{Michael~Tait}
\email{mtait@cmu.edu}
\address{Department of Mathematical Sciences, Carnegie Mellon University, Pittsburgh, PA, 15213, USA}
\thanks{Research of the second author is supported in part by NSF grant DMS-1606350.}

\newtheoremstyle{bfnote}%
{}{}%
{\slshape}{}%
{\bfseries}{\bfseries.}%
{ }%
{\thmname{#1}\thmnumber{ #2}\thmnote{ \ep{\normalfont{}#3}}}

\newtheoremstyle{defbfnote}%
{}{}%
{}{}%
{\bfseries}{.}%
{ }%
{\thmname{#1}\thmnumber{ #2}\thmnote{ (#3)}}

\newtheoremstyle{claim}%
{}{}%
{\slshape}{}%
{\itshape}{.}%
{ }%
{\thmname{#1}\thmnumber{ #2}\thmnote{ \ep{\normalfont{}#3}}}

\theoremstyle{bfnote}
\newtheorem{theo}[equation]{Theorem}

\newtheorem{lemma}[equation]{Lemma}
\newtheorem{corl}[equation]{Corollary}

\newtheorem*{claim*}{Claim}
\newtheorem*{corl*}{Corollary}

\newcommand*{\myproofname}{Proof}

\theoremstyle{definition}

\newtheorem*{defn*}{Definition}

\newtheorem{ques}[equation]{Question}

\newtheorem*{exmp*}{Example}

\newtheorem{prob}[equation]{Problem}

\theoremstyle{remark}
\newtheorem*{ques*}{Question}
\newtheorem*{remk*}{Remark}

\newcommand{\set}[1]{\{#1\}}

\newcommand{\N}{{\mathbb{N}}}
\newcommand{\Z}{\mathbb{Z}}
\newcommand{\F}{\mathbb{F}}
\newcommand{\R}{\mathbb{R}}
\renewcommand{\C}{\mathbb{C}}
\newcommand{\Q}{\mathbb{Q}}

\renewcommand{\epsilon}{\varepsilon}
\renewcommand{\phi}{\varphi}
\renewcommand{\theta}{\vartheta}
\newcommand{\uni}{\mathsf{uni}}
\renewcommand{\Re}{\mathrm{Re}}
\renewcommand{\leq}{\leqslant}
\renewcommand{\geq}{\geqslant}

\renewcommand{\G}{\Gamma}
\renewcommand{\Prob}{\operatorname{Prob}}
\newcommand{\defeq}{\coloneqq}

\newcommand{\emphd}[1]{{\fontseries{b}\selectfont\textsf{#1}}}
\newcommand{\db}{\mathsf{d}}
\newcommand{\Db}{\mathsf{D}}
\newcommand{\T}{\mathbb{T}}
\newcommand{\Diff}{\mathrm{d}}
\newcommand{\bemph}[1]{{\normalfont#1}} 
\newcommand{\ep}[1]{\bemph{(}#1\bemph{)}} 

\newenvironment{scproof}[1][Proof]{\begin{proof}[\textsc{#1}]}{\end{proof}}

\numberwithin{equation}{section}

\renewcommand{\thesubsection}{\arabic{section}.\Alph{subsection}}

\usepackage{etoolbox}

\titleformat{\section}[block]{\scshape\filcenter}{\thesection.}{1ex}{}
\titleformat{\subsection}[block]{\bfseries\filcenter}{\thesubsection.}{1ex}{}
\titleformat{\subsubsection}[runin]{\bfseries}{\thesubsubsection.}{1ex}{}[.]

\titlespacing*{\section}{0pt}{*3}{*1}
\titlespacing*{\subsection}{0pt}{*2}{*1}

\makeatletter
\newcommand{\neutralize}[1]{\expandafter\let\csname c@#1\endcsname\count@}
\makeatother

\renewcommand{\mkbibnamefirst}[1]{\textsc{#1}}
\renewcommand{\mkbibnamelast}[1]{\textsc{#1}}
\renewcommand{\mkbibnameprefix}[1]{\textsc{#1}}
\renewcommand{\mkbibnameaffix}[1]{\textsc{#1}}
\renewcommand{\finentrypunct}{}

\renewbibmacro{in:}{}

\renewbibmacro*{volume+number+eid}{%
	\printfield{volume}%
	\setunit*{\addnbspace}
	\printfield{number}%
	\setunit{\addcomma\space}%
	\printfield{eid}}

\DeclareFieldFormat[article]{volume}{\textbf{#1}\space}
\DeclareFieldFormat[article]{number}{\mkbibparens{#1}}

\DeclareFieldFormat{journaltitle}{#1,}
\DeclareFieldFormat[thesis]{title}{\mkbibemph{#1}\addperiod}
\DeclareFieldFormat[article, unpublished, thesis]{title}{\mkbibemph{#1},}
\DeclareFieldFormat[book]{title}{\mkbibemph{#1}\addperiod}
\DeclareFieldFormat[unpublished]{howpublished}{#1, }

\DeclareFieldFormat{pages}{#1}

\DeclareFieldFormat[article]{series}{Ser.~#1\addcomma}

\begin{document}
	\pagestyle{plain}
	
	\maketitle
	
	\begin{abstract}
		A difference basis with respect to $n$ is a subset $A \subseteq \Z$ such that $A - A \supseteq \set{1, \ldots, n}$. R\'edei and R\'enyi showed that the minimum size of a difference basis with respect to $n$ is $(c+o(1))\sqrt{n}$ for some positive constant $c$. The best previously known lower bound on $c$ is $c \geq 1.5602\ldots$, which was obtained by Leech using a version of an earlier argument due to R\'edei and R\'enyi. In this note we use Fourier-analytic tools to show that the Leech--R\'edei--R\'enyi lower bound is not sharp.
	\end{abstract}
	
	\section{Introduction}
	
	We use $\N$ \ep{resp.\ $\N^+$} to denote the set of all nonnegative \ep{resp.\ positive} integers. For $n \in \N^+$, let $[n] \defeq \set{1, \ldots, n}$ and $[-n] \defeq \set{-n, \ldots, -1}$. Given $A \subseteq \Z$, we write $A - A \defeq \set{a - b \,:\, a,\, b \in A}$.
	
	A set $A \subseteq \Z$ is called a \emphd{difference basis with respect to $n$} if $A - A \supseteq [n]$. In this note we address the following problem, first raised by R\'edei and R\'enyi \cite{Renyi}:
	
	\begin{prob}\label{prob:question}
		For given $n \in \N^+$, what is the minimum size of a difference basis with respect to $n$?
	\end{prob}
	
	Problem \ref{prob:question}, while it is a natural combinatorial number theory question in its own right, also has applications to graceful labelings of graphs \cite{Golomb, GrahamSloane}, to symmetric intersecting families of sets \cite{Bhargav}, and to signal processing \cite{Haykin, Linebarger, Moffet}.
	
	Let $\Db(n)$ denote the smallest size of a difference basis with respect to $n$. In their seminal paper \cite{Renyi}, R\'edei and R\'enyi showed that the limit
	\[
		\db^\ast \,\defeq\, \lim_{n \to \infty} \frac{\Db(n)^2}{n}
	\]
	exists. Clearly, if $[n] \subseteq A - A$, then $n \leq {|A| \choose 2}$, and hence $\db^\ast \geq 2$. On the other hand, it is not hard to give a construction that shows $\db^\ast \leq 4$. It turns out that both these bounds can be improved. In particular, R\'edei and R\'enyi \cite{Renyi} showed that
	\[
		2.4244\ldots \,=\, 2 + \frac{4}{3\pi} \,\leq\, \db^\ast \,\leq\, \frac{8}{3} \,=\, 2.6666\ldots.
	\]
	Leech \cite{Leech} found a way to improve the R\'edei--R\'enyi construction to derive the upper bound $\db^\ast \leq 2.6646\ldots$. This was further improved by Golay \cite{Golay} to $\db^\ast \leq 2.6458\ldots$.
	
	In this note we are interested in lower bounds on $\db^\ast$. Here, again, the result of R\'edei and R\'enyi was improved by Leech \cite{Leech}, who noticed that the argument from \cite{Renyi} depends on a certain parameter $\theta$ \ep{taken by R\'edei and R\'enyi to be $\theta = 3\pi/2$} and that making the optimal choice for $\theta$ gives the following:
	
	\begin{theo}[{Leech--R\'edei--R\'enyi \cite{Leech}}]\label{theo:LRR}
		We have
		\[
			\db^\ast \,\geq\, 2 - 2\inf_{\theta \,\neq\, 0}  \frac{\sin(\theta)}{\theta} \,=\, 2.4344\ldots.
		\]
	\end{theo} 
	
	The contribution of this paper is to show that the bound in Theorem~\ref{theo:LRR} is not sharp:
	
	\begin{theo}\label{theo:main_existential}
		There exists $\epsilon > 0$ such that
		\[
		\db^\ast \,\geq\, \epsilon + 2 - 2\inf_{\theta \,\neq\, 0}  \frac{\sin(\theta)}{\theta}.
		\]
	\end{theo}
	
	Our numerical computations suggest that $\epsilon$ in Theorem~\ref{theo:main_existential} can be taken to be around $10^{-3}$. However, we did not make an effort to optimize $\epsilon$, since it is unclear how close the best lower bound that our methods can give is to the correct value of $\db^\ast$.
	
	Our proof techniques are Fourier-analytic. The original approach of R\'edei and R\'enyi can be formulated in terms of looking at the first Fourier coefficient of a certain probability measure on the unit circle. Essentially, we show that taking into account higher Fourier coefficients leads to better lower bounds on $\db^\ast$.
	
	\section{Preliminaries}
	
	\subsubsection*{Measures}
	
	For a nonempty finite set $A$, $\uni(A)$ denotes the uniform probability measure on $A$. For a function $\phi \colon X \to Y$ and a measure $\mu$ on $X$, the pushforward of $\mu$ by $\phi$ is denoted by $\phi_\ast(\mu)$.
	
	\subsubsection*{The space of measures}
	
	Let $X$ be a compact metric space. We use $\Prob(X)$ to denote the space of all probability Borel measures on $X$ equipped with the usual weak-$\ast$ topology \ep{see, e.g., \cite[\S{}17.E]{K_DST}}. Note that the space $\Prob(X)$ is compact and metrizable \cite[Theorem 17.22]{K_DST}.
	
	\subsubsection*{Measures on the unit circle}
	
	Let $\T \defeq \set{z \in \C \,:\, |z| = 1}$ be the unit circle in the complex plane, viewed as a compact Abelian group. Given a measure $\mu \in \Prob(\T)$, we use $\overline{\mu}$ to denote the pushforward of $\mu$ by the conjugation map $\T \to \T \colon z \mapsto \overline{z}$. The \emphd{Fourier transform} of a measure $\mu \in \Prob(\T)$ is the function $\widehat{\mu} \colon \Z \to \C$ defined by the formula
	\[
	\widehat{\mu}(k) \,\defeq\, \int_{\T} z^k \,\Diff \mu(z).
	\]
	The values $\widehat{\mu}(k)$ are referred to as the \emphd{Fourier coefficients} of $\mu$. We shall make use of the following basic observation:
	
	\begin{lemma}\label{lemma:pos}
		Let $\mu$ be a probability measure on $\T$ and let $A$ be the $n$-by-$n$ matrix with entries \[A(i, j) \defeq \widehat{\mu}(j - i), \qquad \text{for all } 1 \leq i,\, j \leq n.\] Then $A$ is Hermitian and positive semidefinite.
	\end{lemma}
	\begin{scproof}
		That $A$ is Hermitian is clear. To show that $A$ is positive semidefinite, take any $w \in \C^n$. Viewing $w$ as a column vector, we compute
		\begin{align*}
			\left\langle Aw, w \right\rangle \,=\, \sum_{i = 1}^n \sum_{j = 1}^n A(i, j) \overline{w_i} w_j \,&=\, \sum_{i = 1}^n \sum_{j = 1}^n \widehat{\mu}(j - i) \overline{w_i} w_j \,=\, \sum_{i = 1}^n \sum_{j = 1}^n \int_{\T} z^{j-i} \,\Diff \mu(z)  \overline{w_i} w_j \\
			&=\, \int_{\T} \sum_{i = 1}^n \sum_{j = 1}^n \overline{(w_i z^i)} (w_j z^j) \,\Diff \mu(z) \,=\, \int_{\T} \left|\sum_{i=1}^n w_i z^i\right|^2 \,\Diff \mu(z) \,\geq \, 0. \qedhere
		\end{align*}
	\end{scproof}
	
	It will be useful to remember that if a Hermitian matrix $A$ is positive\-/semidefinite, then so is the real symmetric matrix whose entries are the real parts of the corresponding entries of $A$.
	
	For completeness, we record here the converse of Lemma~\ref{lemma:pos} \ep{although we will not need it}:
	
	\begin{theo}[{Bochner--Herglotz \cite[\S1.4.3]{Rudin}}]
		Let $f \colon \Z \to \C$ be a function such that:
		\begin{itemize}
			\item $f(0) = 1$,
			\item $f(-k) = \overline{f(k)}$ for all $k \in \Z$, and
			\item for each $n \in \N^+$, the $n$-by-$n$ matrix $A$ with entries $A(i,j) \defeq f(j-i)$ is positive semidefinite.
		\end{itemize}
	Then there exists a unique probability measure $\mu \in \Prob(\T)$ with $f = \widehat{\mu}$.
	\end{theo}
	
	\subsubsection*{Convolutions of measures} Given two probability measures $\mu$, $\nu$ on $\T$, their \emphd{convolution} is the probability measure $\mu \ast \nu$ on $\T$ given by
	\[
		\int_{\T} f(z) \,\Diff (\mu \ast \nu)(z) \,\defeq\, \int_{\T \times \T} f(xy) \,\Diff (\mu \times \nu)(x, y) \,=\, \int_{\T} \int_{\T} f(xy) \,\Diff\mu(x) \,\Diff \nu(y).
	\]
	Notice that the Fourier transform turns convolution into multiplication, in the sense that \[\widehat{\mu \ast \nu}(k) \,=\, \widehat{\mu}(k) \widehat{\nu}(k) \qquad \text{for all } k \in \Z.\]
	
	\section{Proof of Theorem~\ref{theo:main_existential}}\label{sec:proof}
	
	In this section we prove Theorem~\ref{theo:main_existential}, without making any attempt to compute an exact value for $\epsilon$. Let $\theta = 4.4934\ldots$ be the value for which $\sin(\theta)/\theta$ is minimized \ep{so $\sin(\theta)/\theta = -0.2172\ldots$}. Suppose, towards a contradiction, that there is an infinite set of ``bad'' integers $B \subseteq \N^+$ and a way to assign to every $n \in B$ a difference basis $A_n \subset \Z$ with respect to $n$ so that
	\begin{equation}\label{eq:size}
		|A_n|^2 \,\leq\, \left(2 - \frac{2\sin(\theta)}{\theta} + o(1)\right)n \,=\, (2.4344\ldots + o(1))n.
	\end{equation}
	Take any $n \in B$ and let $\alpha_n \defeq |A_n|^2/n-2$, so $|A_n|^2 = (2 + \alpha_n) n$. Let $\phi_n \colon \Z\to \T$ be the function given by $\phi_n(k) \defeq \exp\left(\theta ik/n\right)$, and define the following two measures on $\T$:
	\[
		\mu_n \defeq (\phi_n)_\ast(\uni(A_n)) \qquad \text{and} \qquad \nu_n \defeq (\phi_n)_\ast(\uni([-n] \cup [n])).
	\]
	
	\begin{lemma}
		For each $n \in B$, there exists a probability measure $\zeta_n\in \Prob(\T)$ such that
		\begin{equation}\label{eq:conv}
			\mu_n \ast \overline{\mu_n} \,=\, \frac{2}{2+\alpha_n} \nu_n + \frac{\alpha_n}{2 + \alpha_n}\zeta_n,
		\end{equation}
	\end{lemma}
	\begin{scproof}
		Let $\xi_n$ be the probability measure on the \ep{finite} set $A_n - A_n$ given by
		\[
			\xi_n(\set{c}) \defeq \frac{1}{|A_n|^2} |\set{(a,b) \in A_n \times A_n \,:\, a-b = c}|.
		\]
		Note that $A_n - A_n \supseteq [-n] \cup [n]$, and hence for each $k \in [-n] \cup [n]$, we have
		\[
			\xi_n(\set{k}) \,\geq\, \frac{1}{|A_n|^2} \,=\, \frac{1}{(2+\alpha_n)n} \,=\, \frac{2}{2+\alpha_n} (\uni([-n] \cup [n]))(\set{k}).
		\]
		It remains to observe that $\mu_n \ast \overline{\mu_n} = (\phi_n)_\ast(\xi_n)$, as
		\begin{align*}
			\int_\T f(z) \,\Diff(\mu_n \ast \overline{\mu_n})(z) \,&=\, \int_{\T\times \T} f(xy) \,\Diff (\mu_n \times \overline{\mu_n})(x, y) \\
			&=\,\frac{1}{|A_n|^2} \sum_{(a,b) \,\in\, A_n \times A_n} f(\phi_n(a) \overline{\phi_n(b)}) \\
			&=\, \frac{1}{|A_n|^2} \sum_{(a,b) \,\in\, A_n \times A_n} f(\phi_n(a-b)) \,=\, \int_\T f(z) \,\Diff (\phi_\ast(\xi_n))(z). \qedhere
		\end{align*}
	\end{scproof}
	
	Now we pass to the limit as $n$ tends to infinity. Let $\phi \colon [-1; 1] \to \T$ be given by $\phi(a) \defeq \exp(\theta i a)$, and let $\nu \defeq \phi_\ast(\lambda)$, where $\lambda$ is the uniform probability measure on $[-1;1]$. It is then clear that \[\nu = \lim_{n \in B} \nu_n.\] Upon replacing $B$ by a subset if necessary, we may also assume that the following limits exist:
	\[
		\alpha \defeq \lim_{n \in B} \alpha_n, \qquad \mu \defeq \lim_{n \in B} \mu_n, \qquad \text{and} \qquad \zeta \defeq \lim_{n \in B} \zeta_n.
	\]
	By \eqref{eq:size}, we have $\alpha \leq -2\sin(\theta)/\theta = 0.4344\ldots$, while from \eqref{eq:conv}, we conclude that
	\begin{equation}\label{eq:conv_lim}
		\mu \ast \overline{\mu} \,=\, \frac{2}{2+\alpha} \nu + \frac{\alpha}{2 + \alpha}\zeta.
	\end{equation}
	
	\begin{lemma}\label{lemma:nu_Fourier}
		The Fourier coefficients of $\nu$ are $\widehat{\nu}(0)=1$ and $\widehat{\nu}(k) = \sin(k\theta)/(k\theta)$ for all $k \neq 0$.
	\end{lemma}
	\begin{scproof}
		A straightforward direct computation.
	\end{scproof}
	
	Let $\delta_1$ denote the Dirac probability measure concentrated at $1 \in \T.$
	
	\begin{corl}\label{corl:sharp}
		The following statements are valid:
		\[\alpha = -2\sin(\theta)/\theta; \qquad \widehat{\mu}(1) = 0; \qquad \text{and} \qquad \zeta = \delta_1.\]
	\end{corl}
	\begin{scproof}
		From \eqref{eq:conv_lim} and Lemma~\ref{lemma:nu_Fourier}, we obtain
		\begin{align}
			0\,\leq\, |\widehat{\mu}(1)|^2 \,=\, \widehat{\mu \ast \overline{\mu}}(1) \,&=\, \frac{2}{2+\alpha} \widehat{\nu}(1) + \frac{\alpha}{2 + \alpha}\widehat{\zeta}(1)\nonumber\\
			&=\, \frac{2}{2+\alpha} \cdot \frac{\sin(\theta)}{\theta} + \frac{\alpha}{2 + \alpha}\widehat{\zeta}(1) \,\leq\, \frac{2}{2+\alpha} \cdot \frac{\sin(\theta)}{\theta} + \frac{\alpha}{2 + \alpha},\label{eq:bound_on_alpha}
		\end{align}
		and therefore $\alpha \geq -2\sin(\theta)/\theta$ \ep{this is essentially the Leech--R\'edei--R\'enyi's proof of Theorem~\ref{theo:LRR}}. Since $\alpha \leq-2\sin(\theta)/\theta$ by assumption, we conclude that $\alpha = -2\sin(\theta)/\theta$ and neither of the two inequalities in \eqref{eq:bound_on_alpha} can be strict, which means that \[\widehat{\mu}(1) = 0 \qquad \text{and} \qquad \widehat{\zeta}(1) = 1.\] Since $\delta_1$ is the only probability measure on $\T$ whose first Fourier coefficient is $1$, we have $\zeta = \delta_1$.
	\end{scproof}
	
	Set $\beta \defeq \sqrt{\alpha/(2+\alpha)} = 0.4224\ldots$. Using Corollary~\ref{corl:sharp}, we can rewrite \eqref{eq:conv_lim} as
	\begin{equation}\label{eq:convv}
		\mu \ast \overline{\mu} \,=\, (1-\beta^2) \nu + \beta^2\delta_1.
	\end{equation}
	
	\begin{lemma}
		The measure $\mu$ has precisely one atom $z \in \T$, and it satisfies $\mu(\set{z}) = \beta$.
	\end{lemma}
	\begin{scproof}
		From \eqref{eq:convv}, it follows that $\mu \ast \overline{\mu}$ has a unique atom, namely $1$, and $(\mu \ast \overline{\mu})(\set{1}) = \beta^2$. If $\mu$ were atomless, then so would be $\mu \ast \overline{\mu}$, so $\mu$ must have at least one atom. On the other hand, if $\mu$ had two distinct atoms, say $x$ and $y$, then we would have $(\mu \ast \overline{\mu})(\set{xy^{-1}}) \geq \mu(\set{x}) \mu(\set{y}) > 0$, which is impossible as $xy^{-1} \neq 1$. Therefore, $\mu$ has a unique atom $z$, and furthermore
		\[
			\mu(\set{z})^2 \,=\, (\mu \ast \overline{\mu})(\set{1}) \,=\, \beta^2,
		\]
		i.e., $\mu(\set{z}) = \beta$, as desired.
	\end{scproof}
	
	If necessary, we may rotate $\mu$ so that its unique atom is $1 \in \T$. Then $\mu$ can be decomposed as
	\begin{equation}\label{eq:eta}
		\mu \,=\, (1 - \beta) \eta + \beta \delta_1,
	\end{equation}
	for some $\eta \in \Prob(\T)$. From \eqref{eq:eta}, we obtain
	\[
		\mu \ast \overline{\mu} \,=\, (1-\beta)^2 (\eta \ast \overline{\eta}) + (1-\beta)\beta(\eta + \overline{\eta}) + \beta^2 \delta_1.
	\]
	Combined with \eqref{eq:convv}, this yields
	\begin{equation}\label{eq:eta1}
		(1-\beta)(\eta \ast \overline{\eta}) + \beta (\eta + \overline{\eta}) = (1+\beta)\nu.
	\end{equation}
	
	\begin{lemma}
		We have $\hat{\eta}(0) = 1$ and $\hat{\eta}(1) = -\beta/(1-\beta) = -0.7314\ldots$.
	\end{lemma}
	\begin{scproof}
		We have $\widehat{\eta}(0) = 1$ since $\eta$ is a probability measure. From \eqref{eq:eta} and Corollary~\ref{corl:sharp}, 
		\[
		0 \,=\, \widehat{\mu}(1) \,=\, (1-\beta) \widehat{\eta}(1) + \beta,
		\]
		which yields $\widehat{\eta}(1) = -\beta/(1 - \beta)$, as desired.
	\end{scproof}
	
	For brevity, set $\gamma \defeq - \beta/(1-\beta)$.
	
	\begin{lemma}\label{lemma:x}
		We have $0 < \Re(\hat{\eta}(2)) < 0.1$.
	\end{lemma}
	\begin{scproof}
		From \eqref{eq:eta1} and Lemma~\ref{lemma:nu_Fourier}, we obtain
		\[
			(1-\beta) |\hat{\eta}(2)|^2 + 2\beta \Re(\hat{\eta}(2)) - (1+\beta) \frac{\sin(2\theta)}{2\theta} \,=\, 0.
		\]
		Setting $x \defeq \Re(\hat{\eta}(2))$, we conclude that
		\[
			(1-\beta)x^2 + 2\beta x - (1+\beta) \frac{\sin(2\theta)}{2\theta} \,\leq\, 0.
		\]
		Using the numerical values for $\beta = 0.4224\ldots$ and $\theta = 4.4934\ldots$, we deduce that
		\[
			-1.5384\ldots \,\leq\, x \,\leq\, 0.0755\ldots \,<\, 0.1.
		\]
		To show that $x > 0$, consider the $3$-by-$3$ matrix $A$ with entries $A(i,j) \defeq \Re(\widehat{\eta}(j - i))$:
		\[
			A \,=\, \left[\begin{array}{ccc}
			1 & \gamma & x \\
			\gamma & 1 & \gamma \\
			x & \gamma & 1
			\end{array}\right].
			\]
		By Lemma~\ref{lemma:pos}, the matrix $A$ must be positive semidefinite. In particular,
		\[
		\det(A) \,=\, (x-1)(-x+2\gamma^2-1) \,\geq\, 0,
		\]
		which yields
		$
			0 < 0.0700\ldots = 2\gamma^2-1 \leq x\leq1
		$.
	\end{scproof}
	
	We are now ready for the final step. Set
	\[
		x \defeq \Re (\hat{\eta}(2)) \qquad \text{and} \qquad y \defeq \Re(\hat{\eta}(3)),
	\]
	and let $M$ be the $4$-by-$4$ matrix with entries $M(i,j) \defeq \Re(\hat{\eta}(j-i))$:
	\[
		M \,=\, \left[\begin{array}{cccc}
		1 & \gamma & x & y \\
		\gamma & 1 & \gamma & x \\
		x & \gamma & 1 & \gamma \\
		y & x & \gamma & 1
		\end{array}\right].
		\]
	By Lemma~\ref{lemma:pos}, the matrix $M$ must be positive semidefinite. In particular,
	\begin{align*}
		\det M \,=\, &\left(\left(-1 - \gamma \right)y + x^2 + 2\gamma x  + \gamma^2 -\gamma-1\right)\\
		&\cdot \, \left(\left(1 - \gamma \right)y + x^2 - 2\gamma x  + \gamma^2 +\gamma-1\right) \,\geq\, 0.
	\end{align*}
	This means that $y$ is located in the interval between
	\[
		y_1 \defeq \frac{x^2 + 2\gamma x + \gamma^2 - \gamma - 1}{\gamma+1} \qquad \text{and} \qquad y_2 \defeq \frac{x^2 - 2\gamma x + \gamma^2 + \gamma - 1}{\gamma-1}. 
	\]
	As a function of $x$, $y_1$ attains its minimum at the point $-\gamma = 0.7314\ldots$. This means that on the interval $[0;0.1]$ it is decreasing, and hence, since $0 < x < 0.1$ by Lemma~\ref{lemma:x}, we conclude that
	\[
		y_1 \,\geq\, \frac{0.01 + 0.2\gamma + \gamma^2 - \gamma - 1}{\gamma+1} \,=\, 0.4848\ldots \,>\, 0.4.
	\]
	Similarly, $y_2$, viewed as a function of $x$, attains its maximum at the point $\gamma = -0.7314\ldots$. Hence, it is decreasing on the interval $[0;0.1]$, and thus
	\[
		y_2 \,\geq\, \frac{0.01 - 0.2\gamma + \gamma^2 + \gamma - 1}{\gamma-1} \,=\, 0.6007\ldots \,>\, 0.4.
	\]
	Therefore, we conclude that $y > 0.4$. On the other hand, from \eqref{eq:eta1} and Lemma~\ref{lemma:nu_Fourier}, we obtain
	\[
	(1-\beta) |\hat{\eta}(3)|^2 + 2\beta \Re(\hat{\eta}(3)) - (1+\beta) \frac{\sin(3\theta)}{3\theta} \,=\, 0,
	\]
	which yields
	\[
		(1-\beta) y^2 + 2\beta y - (1+\beta) \frac{\sin(3\theta)}{3\theta} \,\leq\, 0.
	\]
	Using the numerical values for $\beta = 0.4224\ldots$ and $\theta = 4.4934\ldots$, we obtain
	\[
		-1.5559\ldots \,\leq\, y \,\leq\, 0.0929\ldots \,<\,0.1.
	\]
	This contradiction completes the proof of Theorem~\ref{theo:main_existential}.
	
	\section*{Concluding remarks and acknowledgments}
	
	Even though our proof, as presented in Section~\ref{sec:proof}, does not give an explicit lower bound on $\epsilon$, it is clear how one could obtain such an explicit lower bound by introducing small margins of error throughout the argument. However, determining the optimal value of $\epsilon$ in Theorem \ref{theo:main_existential} appears technically challenging. One difficulty is that is is necessary to quantify how ``close'' the measure $\zeta$ is to the Dirac measure in Corollary \ref{corl:sharp}; the outcome of this step then propagates through the rest of the proof. It seems unlikely that our methods could yield the exact value of $\db^\ast$. Golay felt that the correct value ``will, undoubtedly, never be expressed in closed form'' \cite{Golay}. Nevertheless, we do not know the answer to the following question:
	
	\begin{ques}
		Let $\mathsf{a}$ denote the infimum of all real numbers $\alpha > 0$ such that there exist probability measures $\mu$, $\zeta \in \Prob(\T)$ satisfying \eqref{eq:conv_lim}. We know that $\db^\ast \geq 2 + \mathsf{a}$. Is it true that, in fact, $\db^\ast = 2 + \mathsf{a}$?
	\end{ques}
	
	The second author would like to thank Craig Timmons for introducing him to the problem. We are very grateful to the anonymous referee for carefully reading the manuscript and providing helpful suggestions. 
	
		\printbibliography

\end{document}